\documentclass[12pt]{article} 
\usepackage[final]{graphicx}
\usepackage[dvips]{color}
\usepackage{graphicx}
\usepackage{amsmath}
\usepackage{amssymb}
\usepackage{ascmac}
\usepackage{amsthm}
\usepackage{bm}
\usepackage{fancyhdr}
\def\Bbb{\bf} 

\newcommand\R{{\Bbb R}}

\newcommand\Q{{\Bbb Q}}

\newtheorem{lmm}{Lemma}[section]
\newtheorem{thm}[lmm]{Theorem}

\newtheorem{rmk}[lmm]{Remark}

\newtheorem{qst}[lmm]{Question}
\newtheorem{ans}[lmm]{Answer}
\def\comment#1{ }

\makeatletter
    
    \@addtoreset{equation}{section}
  \makeatother
\allowdisplaybreaks
\title{An approximation of a catenoid \\
constructed from  \\
piecewise truncated conical minimal surfaces}
\author{Akihito Ebisu and Yoshiroh Machigashira
}
\date{\today}
\begin{document}
\maketitle
\begin{abstract}
In [3], we considered an approximation of a catenoid 
constructed from even truncated cones that maintains minimality in a certain sense.
In this paper, 
we consider such an approximation consisting of odd truncated cones that maintains minimality in the same sense. 
Through this procedure, we obtain a discrete curve approximating a catenary by exploiting the fact that it is the function
that generates a catenoid.
In this investigation, the theory of the Gauss hypergeometric functions plays an important role.

\end{abstract}

\section{Introduction}
In this paper, we consider the problem of approximating a catenoid, which is a minimal surface,
using odd truncated cones in such a manner that maintains {\it minimality} in a certain sense. 
From this surface, we also obtain an approximation of a catenary in the form of a polyline
that retains what we regard to be the most important property of the catenrary.

As a real-world illustration of the problem we consider, suppose we wish to form a surface from a rubber membrane.
We can consider such a surface to be ``stable'' if it is difficult to deform, i.e.,
if it resists compression and tension. 
Only minimal surfaces are stable in this sense.
For this reason, the study of minimal surfaces is important for industrial applications.
Now, suppose that instead of a rubber membrane, we wish to construct a surface from plane figures. 
This is a common situation in industrial applications.
A truncated cone can be constructed from plane figures, and thus,
for such applications, 
it would be useful to develop methods for approximating surfaces of various types using truncated cones.
Further, if we could construct surfaces that approximate minimal surfaces
in such a manner that maintains {\it minimality},
they would be of great practical usefulness.
Formulating such a method for constructing approximations of surfaces from truncated cones
would also be important mathematically, because taking the limit of an infinite number of truncated cones,
we could obtain the parametric equations of various surfaces of interest.

In this work,
we approach the problem of approximating a catenoid by considering approximations of the corresponding catenary.
More precisely, we consider polylines forming discrete curves that converge to the catenary.
Further, we choose these polylines in such a way 
that the surfaces of revolution they generate are themselves {\it minimal}.
Because such a polyline represents a discretization of a catenary,
the problem we consider is also useful from a mathematical point of view.
The function $h_{n,2/(2m+1)}(y) (0 \leq n \leq m)$ (cf. Theorem 2.4) used to construct these polylines (cf. (2.2)) has two noteworthy characteristics:
It is a rational function over the rational number field {\Q} for any $n$ and $m$, and it has a closed-form expression.
It is interesting that we can approximate a catenary with polylines obtained using such functions.

\section{Main Theorem}
A catenary is the curve assumed by a hanging chain.
Its form is given by the function $C_c:t \mapsto c\cosh(t/c)$ for a constant $c>0$.
The surface of revolution generated by rotating a catenary $C_c(t)$ about the $t$-axis
is called a catenoid. Here, such a surface is denoted by $R(C_c(t))$.
The catenoid has the following special property:
Every non-planar rotationally symmetric minimal surface
is congruous to a piece of a catenoid (cf. 3.5.1 in [2]). Here, the term
``minimal'' means ``of mean curvature zero."
Now, we consider catenoids $R(C_c|_{(-1,1)})$, generated by catenaries satisfying $t \in [-1,1]$, whose boundaries consist of two circles of radii $a$. 
If we choose $a$ to be sufficiently large, then there are two such catenoids, generated by catenaries whose values of $c$ we write
$c_{a}^{\pm}$, where we choose $c_a^-<c_a^+$.
It is known that the area of $R(C_{c_{a}^{+}}|_{(-1,1)})$ is invariant
with respect to infinitesimal perturbations of the surface that keep the boundaries fixed.
The same holds for $R(C_{c_{a}^{-}}|_{(-1,1)})$
(cf. Theorem 1 in 2.1 of [2]). 
Moreover, it is known that 
the area of $R(C_{c_{a}^{+}}|_{(-1,1)})$ is minimal among the set of surfaces possessing the same boundaries,
while that of $R(C_{c_{a}^{-}}|_{(-1,1)})$ is not. Below, we consider discretizations of $R(C_{c_{a}^{\pm}}|_{(-1,1)})$.

We call a cone whose apex is cut off by a plane parallel to its base a truncated cone.
For $x_0, x_1>0$ and $\ell >0$, let $D_{1,\ell}(x_0,x_1)$ be the truncated cone 
whose circular edges have radii $x_0$ and $x_1$
and whose height is $\ell$.
Here, we do not consider the regions interior to the two circles of radii
$x_0$ and $x_1$ to be part of $D_{1,\ell}(x_0,x_1)$.
Defining 
$$
S_{1,\ell}(x_0,x_1):=(x_0+x_1)\sqrt{(x_1-x_0)^2+\ell ^2},
$$
we see that the area of $D_{1,\ell}(x_0,x_1)$ is equal to $\pi \cdot S_{1,\ell}(x_0,x_1)$.
For $x_0, x_1, x_2>0$ and $\ell >0$, let $D_{2,\ell}(x_0,x_1,x_2)$ be the figure 
consisting of the union of $D_{1,\ell}(x_0,x_1)$ and $D_{1,\ell}(x_1,x_2)$ 
attached along the circle of radius $x_1$.
We define $D_{n,\ell}(x_0,x_1,$ $x_2,\cdots ,x_n)$ for $n \geq 3$ similarly.
The surface $D_{n,\ell}(x_0,x_1,\cdots ,x_n)$ is called a {\it piecewise truncated conical surface} with length $(n;\ell)$,
or simply a {\rm PTC} surface with {\rm L}-$(n;\ell)$.
Note that the boundary  of $D_{n,\ell}(x_0,x_1,\cdots ,x_n)$ consists of two circles of radii $x_0$ and $x_n$,
and its area is given by
$$
\pi \sum _{i=1} ^{n}S_{1,\ell}(x_{i-1},x_{i}).
$$
Next, we define 
$$
S_{n,\ell}(x_{0},x_{1},\cdots x_{n}):=\sum _{i=1} ^{n}S_{1,\ell}(x_{i-1},x_{i}).
$$
For arbitrary fixed $a, b >0$ and $n \in {\Bbb N}:=\{0,1,\cdots \}$,
$D_{n+2,\ell}(a,x_0,x_1,\cdots ,x_n,b)$ is called a {\rm PTC} surface 
with boundary condition $(a,b)$ and length $(n+2, \ell)$, or simply {\rm BCL}-$(a,b;n+2;\ell)$.
A {\rm PTC} surface $D_{n+2,\ell}(a,x_0^{(0)},x_1^{(0)},\cdots ,x_n^{(0)},b)$ with {\rm BCL}-$(a,b;n+2;\ell)$
is said to be {\it minimal} if $(x_0^{(0)},x_1^{(0)},\cdots ,x_n^{(0)})$ is a critical point of the function
$$
(x_0,x_1,\cdots ,x_n)\mapsto S_{n+2,\ell}(a,x_{0},x_{1},\cdots x_{n},b).
$$ 
Moreover, a {\rm PTC} minimal surface $D_{n+2,\ell}(a,x_0^{(0)},x_1^{(0)},\cdots ,x_n^{(0)},b)$
with {\rm BCL}-$(a,b;n+2;\ell)$ is said to be {\it stable} 
if and only if 
the Hessian matrix of the above function is positive definite at $(x_0^{(0)},x_1^{(0)},\cdots ,x_n^{(0)})$.

Studying a symmetric {\rm PTC} surface $D_{2m,\ell}(a,x_{m-1},x_{m-2},\cdots ,x_0,\cdots ,x_{m-2}, x_{m-1},$ $ a)$
with {\rm BCL}-$(a,a;2m;\ell)$, the following questions come to mind:
\begin{qst} \mbox{} \\
1. Where does this surface become minimal? \\
2. Where does this surface become stable?
\end{qst}
We answered these questions in [3].
The answer to the first is the following:
\begin{ans}(cf. Theorem 1 in [3])\\
For $n \in {\Bbb N}$ and $\ell > 0$, there exist an explicit function $g_{n,\ell}(x)$ on $\R_{>0}$ 
and $\eta_{n,\ell}$
such that for any $m \in {\Bbb N}$, the following holds:
\begin{description}
\item[(1)] If $a> \eta _{m,\ell}$, then the equation $g_{m,\ell}(x)-a=0$ has two positive solutions
$x_{a,m,\ell}^{\pm}$ with $x_{a,m,\ell}^{-} < x_{a,m,\ell}^{+}$.
\item[(2)]$D_{2m,\ell}(a,x_{m-1},x_{m-2},\cdots ,x_0,\cdots ,x_{m-2}, x_{m-1},$ $ a)$
with $x_n=g_{n,\ell}(x_{a,m,\ell}^{\pm})$ for $n=0, \cdots ,m-1$
are {\rm PTC} minimal surfaces.
\end{description}
Moreover, we have
\begin{align*}
g_{n,\ell}(x)&=xT_{n}\left(1+\frac{\ell ^2}{2x^2}\right)=xF\left(n,-n;\frac{1}{2};\frac{-\ell ^2}{4x^2}\right) \\
&=\frac{x}{2}\left\{ \left( \sqrt{1+\left(\frac{\ell}{2x}\right)^2}+\frac{\ell}{2x}\right)^{2n}
+\left( \sqrt{1+\left(\frac{\ell}{2x}\right)^2}-\frac{\ell}{2x}\right)^{2n} 
\right\} ,
\end{align*}
where $T_n$ is a Chebyshev polynomial of the first kind, and $F$ is the Gauss hypergeometric series.
\end{ans}
The answer to the second question is the following:
\begin{ans}(cf. Theorem 2 in [3])
Under the same conditions as in Answer 2.2, $D_{2m,\ell}(a,x_{m-1},x_{m-2},\cdots ,x_0,\cdots ,x_{m-2}, x_{m-1},$ $ a)$ with $x_n=g_{n,\ell}(x_{a,m,\ell}^{+})$ for $n=0, \cdots ,m-1$ is stable.
\end{ans}

In this paper, we consider the same questions 
for another symmetric {\rm PTC} surface, namely
$D_{2m+1,\ell}(a,y_{m-1},y_{m-2},\cdots ,y_0,y_0,\cdots ,y_{m-2}, y_{m-1}, a)$
with {\rm BCL}-$(a,a;2m+1;\ell)$.
A calculation provides the following answer to the first question:
\begin{thm}
For $n \in {\Bbb N}$ and $\ell > 0$, there exist an explicit function 
$h_{n,\ell}(y)$ on $\{y\in \R_{>0};y > \ell /2\}$ 
and $\nu_{n,\ell}$
such that for any $m \in {\Bbb N}$, the following holds:
\begin{description}
\item[(1)] If $a> \nu _{m,\ell}$, then the equation $h_{m,\ell}(y)-a=0$ has two positive solutions
$y_{a,m,\ell}^{\pm}$ with $y_{a,m,\ell}^{-} < y_{a,m,\ell}^{+}$.
\item[(2)] $D_{2m+1,\ell}(a,y_{m-1},y_{m-2},\cdots ,y_0,y_0,\cdots ,y_{m-2}, y_{m-1}, a)$
with $y_n=h_{n,\ell}(y_{a,m,\ell}^{\pm})$ for $n=0, \cdots ,m-1$
are {\rm PTC} minimal surfaces.
\end{description}
Moreover, we have
\begin{align*}
h_{n,\ell}(y)&=yV_{n}\left(1+\frac{2\ell ^2}{4y^2-\ell ^2}\right)
=yF\left(n+1,-n;\frac{1}{2};\frac{-\ell ^2}{4y^2-\ell ^2}\right)\\
&=\frac{(2y+\ell)^{2n+1}+(2y-\ell)^{2n+1}}{4(4y^2-\ell ^2)^n},
\end{align*}
where $V_n$ is a Chebyshev polynomial of the third kind.
\end{thm}
The answer to the second question is
\begin{thm}
Under the same conditions as in Theorem 2.4, 
$D_{2m+1,\ell}(a,y_{m-1},y_{m-2},$ $\cdots ,y_0,y_0,\cdots ,y_{m-2}, y_{m-1}, a)$
with $y_n=h_{n,\ell}(y_{a,m,\ell}^{+})$ for $n=0, \cdots ,m-1$
is stable.
\end{thm}
Now, the {\rm PTC} minimal surfaces considered in Answer 2.2 with $\ell=1/m$ are denoted by $E( x_{a,m,1/m}^{\pm} )$,
and those considered in Theorem 2.4 with $\ell =2/(2m+1)$ are denoted by $E(y_{a,m,2/(2m+1)}^{\pm})$.
From the fact that the areas of $R(C_{c_a^{+}}|_{(-1,1)})$ and $R(C_{c_a^{-}}|_{(-1,1)})$ are both invariant
with respect to infinitesimal perturbation,  
Answers 2.2 and 2.3 and  Theorems 2.4 and 2.5,
we find that as $m\rightarrow \infty$,
$E( x_{a,m,1/m}^{+} )$ and $E(y_{a,m,2/(2m+1)}^{+})$ tend to $R(C_{c_a^{+}}|_{(-1,1)})$,
while $E( x_{a,m,1/m}^{-} )$ and $E(y_{a,m,2/(2m+1)}^{-})$ tend to $R(C_{c_a^{-}}|_{(-1,1)})$.
Therefore, we can regard the polylines whose vertices are specified by the sequence
\begin{gather}
\begin{split}
\left(
(-1, a), \left(\frac{1-m}{m}, g_{m-1,\frac{1}{m}}\left( x_{a,m,\frac{1}{m}}^{\pm} \right)\right) , \cdots ,
\left(0, x_{a,m,\frac{1}{m}}^{\pm}\right) ,\right. \\
\left. \cdots ,
\left(\frac{m-1}{m}, g_{m-1,\frac{1}{m}}\left( x_{a,m,\frac{1}{m}}^{\pm} \right)\right)
,(1,a)
\right)
\end{split}
\end{gather}
and 
\begin{gather}
\begin{split}
\left(
(-1, a), \left(\frac{1-2m}{2m+1}, h_{m-1,\frac{2}{2m+1}}\left( y_{a,m,\frac{2}{2m+1}}^{\pm} \right)\right) , \cdots ,
\left(\frac{-1}{2m+1}, y_{a,m,\frac{2}{2m+1}}^{\pm}\right) ,\right. \\
\left.
\left(\frac{1}{2m+1}, y_{a,m,\frac{2}{2m+1}}^{\pm}\right),\cdots, 
\left(\frac{2m-1}{2m+1}, h_{m-1,\frac{2}{2m+1}}\left( y_{a,m,\frac{2}{2m+1}}^{\pm} \right)\right),(1,a)
\right)
\end{split}
\end{gather}
as {\it good} approximations of the catenaries on $[-1, 1]$ with boundary $a$.

We prove Theorem 2.4 in \S 3 and Theorem 2.5 in \S 4.
In \S 5, by investigating specific catenaries,
we determine the precision of the approximations of the catenary provided by the above polylines.

\begin{rmk}
We showed that $E( x_{a,m,1/m}^{\pm} )$ and $E(y_{a,m,2/(2m+1)}^{\pm})$, respectively, tend to a catenoid
by using the fact that the areas of $R(C_{c_a^{\pm}}|_{(-1,1)})$ are both invariant
with respect to infinitesimal perturbation.
However, we can demonstrate the same thing without using this fact,
since we have obtained
$g_{n,\frac{1}{m}}(x_{a,m,1/m}^{\pm})$ and \\$h_{n,\frac{2}{2m+1}}$ $(y_{a,m,2/(2m+1)}^{\pm})$ explicitly.
Indeed, writing $t:=n/m$, we have
\begin{gather*}
\lim _{m\rightarrow \infty}g_{n,\frac{1}{m}}(x)
=x\cosh (t/x)=C_x(t).
\end{gather*}
Moreover, we know that $x_{a,m,1/m}^{+}$ and $x_{a,m,1/m}^{-}$ tend to $c_a^{+}$ and $c_a^{-}$ as $m\rightarrow \infty$,
respectively, because we have
\begin{gather*}
\lim _{m\rightarrow \infty}g_{m,\frac{1}{m}}(x)
=x\cosh (1/x)=C_x(1).
\end{gather*}
In this way, by constructing PTC minimal surfaces, we can obtain the parametric equations of the corresponding minimal surfaces.
The same is true in the case of $h_{n,\frac{2}{2m+1}}(y_{a,m,2/(2m+1)}^{\pm})$.
\end{rmk}
\section{A proof of Theorem 2.4}
We denote $S_{1,\ell}(s,t)$ by $S_{\ell}(s,t)$ for simplicity.
Recall that
$$
S_{\ell}(s,t)=(s+t)\sqrt{(t-s)^2+\ell ^2}.
$$
Putting 
$$
T_{a,m,\ell}(y_0,y_1,\cdots ,y_{m-1}):=y_0\ell +\sum _{i=1}^{m-1}S_{\ell}(y_{i-1},y_{i})+S_{\ell}(y_{m-1},a),
$$
The area of $D_{2m+1,\ell}(a,y_{m-1},y_{m-2},\cdots ,y_0,y_0,\cdots ,y_{m-2}, y_{m-1}, a)$
is expressed as $2\pi $ $T_{a,m,\ell}(y_0,y_1,\cdots ,y_{m-1})$.
Therefore, we only have to evaluate $T_{a,m,\ell}(y_0,y_1,\cdots ,y_{m-1})$
for investigating minimality and stability of this {\rm PTC} surface.

In this section, we prove Theorem 2.4.

\subsection{The case $m=1$}
For $a>0$ and $\ell >0$, we consider the critical points of the function
$T_{a,1,\ell}(y_0)\left(=y_0\ell +\right.$ $\left. S_{\ell}(y_0,a)\right)$.
Note that
\begin{gather*}
\frac{\partial S_\ell}{\partial s}(s,t)=\frac{2s^2-2ts+\ell ^2}{\sqrt{(t-s)^2+\ell ^2}},\ 
\frac{\partial S_\ell}{\partial t}(s,t)=\frac{2t^2-2st+\ell ^2}{\sqrt{(t-s)^2+\ell ^2}}.
\end{gather*}
Thus, if $y_0$ is a critical point of $T_{a,1,\ell}(y_0)$, then
\begin{gather}
\frac{dT_{a,1,\ell}}{dy_0}(y_0)=\ell +\frac{2y_0^2-2ay_0+\ell ^2}{\sqrt{(a-y_0)^2+\ell ^2}}=0.
\end{gather}
Therefore, we have
$$
(y_0-a)(4y_0^3-4ay_0^2+3\ell ^2y_0 +a\ell ^2)=0.
$$
Because $y_0=a$ does not satisfy (3.1), we obtain
\begin{gather}
a=y_0+\frac{4\ell ^2y_0}{4y_0^2-\ell ^2}.
\end{gather}
If $\ell /2 >y_0>0$ in (3.2), then $a$ is a negative number.
Therefore, we regard the right-hand side of (3.2) as a function on $\{y \in \R_{>0};y>\ell /2\}$,
and denote this function by $h_{1,\ell}(y)$, that is, 
$$
h_{1,\ell}(y):=y+\frac{4\ell ^2y}{4y^2-\ell ^2}.
$$
Since the second derivative of $h_{1,\ell}(y)$ is
$$
\frac{32\ell ^2 y (4y^2+3\ell ^2)}{(4y^2-\ell ^2)^3},
$$
this funtion is positive and convex on $\{y \in \R_{>0};y>\ell /2\}$.
Moreover, because
$$
\lim_{y\rightarrow \ell /2}h_{1,\ell}(y)=\lim_{y\rightarrow \infty}h_{1,\ell}(y)=\infty,
$$
$h_{1,\ell}(y)$ takes the unique minimal value $\nu _{1,\ell}(>0)$ at a point $\mu _{1,\ell}(>\ell/2)$.
Hence, if $a> \nu_{1,\ell}$, then there are the two solutions $y_{a,1,\ell}^{\pm}$ of $h_{1,\ell}(y)-a=0$,
where $y_{a,1,\ell}^{+}>\mu _{1,\ell}>y_{a,1,\ell}^{-}$.
Therefore, $y_{a,1,\ell}^{\pm}$ are the critical points of the function $T_{a,1,\ell}(y_0)$.
Namely, $D_{3,\ell}(a,y_{a,1,\ell}^{\pm},y_{a,1,\ell}^{\pm},a)$ are minimal.

\subsection{The case $m=2$}
For $a>0$ and $\ell >0$, we consider the critical points of the function
$T_{a,2,\ell}(y_0,y_1)\left(=y_0\ell \right.$ $\left. +\right.$ $\left. S_{\ell}(y_0,y_1)+S_{\ell}(y_1,a)\right)$,
that is, we consider a point $(y_0,y_1)$ satisfying 
$$
\frac{\partial T_{a,2,\ell}}{\partial y_0}(y_0,y_1)=\frac{\partial T_{a,2,\ell}}{\partial y_1}(y_0,y_1)=0.
$$
From the formula 
$$
\frac{\partial T_{a,2,\ell}}{\partial y_0}(y_0,y_1)=0,
$$
we see that $y_1=h_{1,\ell}(y_0)$.
Moreover, 
\begin{gather}
0=\frac{\partial T_{a,2,\ell}}{\partial y_1}(y_0,y_1)
=\frac{2y_1^2-2y_0y_1+\ell ^2}{\sqrt{(y_1-y_0)^2+\ell ^2}}+\frac{2y_1^2-2ay_1+\ell ^2}{\sqrt{(a-y_1)^2+\ell ^2}}
\end{gather}
implies that 
\begin{align*}
0&=(2y_1^2-2y_0y_1+\ell ^2)^2((a-y_1)^2+\ell ^2)-(2y_1^2-2ay_1+\ell ^2)^2((y_1-y_0)^2+\ell ^2)\\
&=\ell ^2(a-y_0)(4y_1^3-4ay_0y_1+2\ell ^2 y_1+\ell ^2 y_0 +a\ell ^2).
\end{align*}
Here, if $a=y_0$, then Formula (3.3) does not hold and so we have
$$
(\ell ^2-4y_0y_1)a+4y_1^3+2\ell ^2 y_1+\ell ^2 y_0=0,
$$
that is, 
\begin{align*}
a&=\frac{4y_1^3+2\ell ^2 y_1+\ell ^2 y_0}{4y_0y_1-\ell ^2}
=\frac{4(h_{1,\ell}(y_0))^3+2\ell ^2 h_{1,\ell}(y_0)+\ell ^2 y_0}{4y_0h_{1,\ell}(y_0)-\ell ^2}\\
&=\frac{16y_0^5+40\ell ^2y_0^3+5\ell ^4 y_0}{(4y_0 ^2-\ell ^2)^2}
=y_0+12y_0\frac{\ell ^2}{4y_0^2-\ell ^2}+16y_0\left(\frac{\ell ^2}{4y_0^2-\ell ^2}\right)^2.
\end{align*}
If we put 
\begin{gather}
h_{2,\ell}(y):=y+12y\frac{\ell ^2}{4y^2-\ell ^2}+16y\left(\frac{\ell ^2}{4y^2-\ell ^2}\right)^2,
\end{gather}
then $h_{2,\ell}(y)$ is positive, convex on $\{y\in \R_{>0};y>\ell /2\}$ and
$$
\lim _{y\rightarrow \ell /2}h_{2,\ell}(y)=\lim _{y\rightarrow \infty}h_{2,\ell}(y)=\infty .
$$
Thus, $h_{2,\ell}(y)$ takes the unique minimal value $\nu _{2,\ell}$ at a point $\mu _{2,\ell}$.
Hence, if $a> \nu_{2,\ell}$, then there are the two solutions $y_{a,2,\ell}^{\pm}$ of $h_{2,\ell}(y)-a=0$,
where $y_{a,2,\ell}^{+}>\mu _{2,\ell}>y_{a,2,\ell}^{-}$.
Consequently, if $a> \nu_{2,\ell}$, there are the two critical points
$(y_{a,2,\ell}^{\pm},h_{1,\ell}(y_{a,2,\ell}^{\pm}))$ of $T_{a,2,\ell}(y_0,y_1)$.
Namely, if $a> \nu_{2,\ell}$, then $D_{5,\ell}(a,h_{1,\ell}(y_{a,2,\ell}^{\pm})), y_{a,2,\ell}^{\pm},
y_{a,2,\ell}^{\pm}, h_{1,\ell}(y_{a,2,\ell}^{\pm})), a)$ are minimal.

\subsection{The case $m=3$}
We consider the critical points of
$$
T_{a,3,\ell}(y_0,y_1,y_2)=y_0\ell+S_\ell(y_0,y_1)+S_\ell(y_1,y_2)+S_\ell(y_2,a)
$$
for $a>0$ and $\ell >0$.
If $(y_0,y_1,y_2)$ is a critical point of $T_{a,3,\ell}$, 
then as in the case $m=2$, we have
$$
y_1=h_{1,\ell}(y_0),\ y_2=h_{2,\ell}(y_0),
$$
and 
\begin{align*}
a&=\frac{4y_2^3+2\ell ^2 y_2+\ell ^2 y_1}{4y_1y_2-\ell ^2}
=\frac{4(h_{2,\ell}(y_0))^3+2\ell ^2 h_{2,\ell}(y_0)+\ell ^2 h_{1,\ell}(y_0)}
{4h_{1,\ell}(y_0)h_{2,\ell}(y_0)-\ell ^2}\\
&=\frac{64y_0^7+336\ell ^2 y_0^5+140\ell ^4 y_0^3+7\ell ^6 y_0}{(4y_0^2-\ell ^2)^3}\\
&=y_0+24y_0\left(\frac{\ell ^2}{4y_0 ^2-\ell ^2}\right)+80y_0\left(\frac{\ell ^2}{4y_0 ^2-\ell ^2}\right)^2
+64y_0\left(\frac{\ell ^2}{4y_0 ^2-\ell ^2}\right)^3.
\end{align*}
Putting
\begin{gather}
h_{3,\ell}(y):=y+24y\left(\frac{\ell ^2}{4y ^2-\ell ^2}\right)+80y\left(\frac{\ell ^2}{4y ^2-\ell ^2}\right)^2
+64y\left(\frac{\ell ^2}{4y ^2-\ell ^2}\right)^3,
\end{gather}
similarly as in the case $m=2$, we see that there is $\mu _{3,\ell}>0$ with $h_{3,\ell}'(\mu _{3,\ell})=0$
such that if $a> \nu _{3,\ell}:=h_{3,\ell}(\mu _{3,\ell})$,
then the equation $h_{3,\ell}(y)=a$ has the two solutions $y_{a,3,\ell}^{\pm}$ with
$y_{a,3,\ell}^{+}>\mu _{3,\ell}>y_{a,3,\ell}^{-}$.
Consequently, if $a>\nu _{3,\ell}$, then 
$(y_{a,3,\ell}^{\pm},h_{1,\ell}(y_{a,3,\ell}^{\pm}),h_{2,\ell}(y_{a,3,\ell}^{\pm}))$
are the critical points of $T_{a,3,\ell}(y_0,y_1,y_2)$.

Repeating the above argument,
we see that $h_{n,\ell}(y)$ is defined as
\begin{gather}
h_{n,\ell}(y):=\frac{4(h_{n-1,\ell}(y))^3+2\ell ^2 h_{n-1,\ell}(y) + \ell ^2 h_{n-2,\ell}(y)}
{4h_{n-2,\ell}(y)h_{n-1,\ell}(y)-\ell ^2}
\end{gather}
for $n\geq 2$ and
\begin{align}
h_{0,\ell}(y)&=y,\\
h_{1,\ell}(y)&=y+4y\frac{\ell ^2}{4y^2-\ell ^2}.
\end{align}
From Formulas (3.4), (3.5), (3.7) and (3.8), we can anticipate that 
$h_{n,\ell}(y)$ is expressed as, for $n\in {\Bbb N}$,
\begin{gather}
h_{n,\ell}(y)=yF\left(n+1,-n;\frac{1}{2};\frac{-\ell ^2}{4y^2-\ell ^2}\right),
\end{gather}
where $F(\alpha,\beta;\gamma;x)$ is called the Gauss hypergeometric series and is defined as
$$
F(\alpha,\beta;\gamma;x):=\sum _{i=0}^{\infty}\frac{(\alpha)_i(\beta)_i}{(\gamma)_i(1)_i}x^i
$$
and $(\alpha)_i:=\Gamma(\alpha +i)/\Gamma(\alpha)$, and so on.
Indeed, we prove this Formula (3.9) in the next subsection.

\subsection{Expression of $h_{n,\ell}(y)$}
In this subsection, we prove that $h_{n,\ell}(y)$ is expressed as (3.9) for $n\in {\Bbb N}$.
It is obvious that (3.9) with $n=0$ and $n=1$ are equal to (3.7) and (3.8), respectively. 
Therefore, we have only to show that Formula (3.9) satisfies (3.6) when $n\geq 2$.

The Jacobi polynomial $P_{n}^{(\alpha ,\beta)}(x)$ is expressed as
$$
P_{n}^{(\alpha ,\beta)}(x)=\frac{(\alpha +1)_n}{(1)_n}F\left(\alpha +\beta +1+n,-n;\alpha +1;\frac{1-x}{2} \right)
$$
in terms of the Gauss hypergeometric series(cf. 15.4.6 in [1]). 
Hence, by this relation, we see that the right hand side of (3.9) is equal to 
\begin{gather}
\frac{(1)_n}{(1/2)_n}yP_{n}^{(-1/2,1/2)}\left( 1+\frac{2\ell ^2}{4y^2 -\ell ^2}\right).
\end{gather}
In addition, by using the third kind Chebyshev polynomial $V_n(x)$ which is expressed as
$$
V_n(x)=\frac{2^{2n}(1)_n(1)_n}{(1)_{2n}}P_n^{(-1/2,1/2)}(x),
$$
(cf. 1.2.3 in [4]), we see that (3.10) can be expressed as
$$
yV_n\left( 1+\frac{2\ell ^2}{4y^2 -\ell ^2}\right).
$$
Consequently, we found that the right hand side of (3.9) is rewritten in terms of the third kind Chebyshev polynomial
and denote this by $\tilde{h}_{n,\ell}(y)$, that is,
\begin{gather}
\tilde{h}_{n,\ell}(y):=yV_n\left( 1+\frac{2\ell ^2}{4y^2 -\ell ^2}\right).
\end{gather}
So, it suffices to prove that $\tilde{h}_{n,\ell}(y)$ satisfies the same formula as (3.6) when $n\geq 2$.

We remark that the third kind Chebyshev polynomial satisfies the following relation:
\begin{gather}
V_n(x)=2xV_{n-1}(x)-V_{n-2}(x),
\end{gather}
where $V_0(x)=1, V_1(x)=2x-1$(cf. 1.2.3(1.12a) in [4]). Then, we see that 
\begin{gather}
V_{n-1}^2(x)-V_n(x)V_{n-2}(x)=2-2x
\end{gather}
for $n\geq 2$. Formula (3.13) is showed in the same way as Lemma 1 in [3].

Rearranging (3.6), we see that the formula we should show is
\begin{gather}
4\tilde{h}_{n-1,\ell}(y)\left( \tilde{h}_{n-1,\ell}^2(y)-\tilde{h}_{n,\ell}(y)\tilde{h}_{n-2,\ell}(y)\right)
+\ell ^2 \left( \tilde{h}_{n,\ell}(y)+2\tilde{h}_{n-1,\ell}(y)+\tilde{h}_{n-2,\ell}(y)\right)=0.
\end{gather}
Substituting (3.11) for each $\tilde{h}_{i,\ell}(y$) in the left hand side of (3.14),
and using (3.12) and (3.13),
we see that the left hand side of (3.14) is equal to zero.
Consequently, we have
$$
h_{n,\ell}(y)=\tilde{h}_{n,\ell}(y)=yV_{n}\left(1+\frac{2\ell ^2}{4y^2-\ell ^2}\right)
=yF\left(n+1,-n;\frac{1}{2};\frac{-\ell ^2}{4y^2-\ell ^2}\right).
$$
Moreover, Because 
$$
F\left(n+1, -n; \frac{1}{2}; -x^2 \right) =
\frac{\left( \sqrt{1+x^2}+x\right) ^{2n+1}+\left( \sqrt{1+x^2}-x\right) ^{2n+1}}{2\sqrt{1+x^2}}
$$
(cf. 15.1.12), we obtain
\begin{gather*}
h_{n,\ell}(y)=yF\left(n+1,-n;\frac{1}{2};\frac{-\ell ^2}{4y^2-\ell ^2}\right)
=\frac{(2y+\ell)^{2n+1}+(2y-\ell)^{2n+1}}{4(4y^2-\ell ^2)^n}.
\end{gather*}
\subsection{A proof of Theorem 2.4}
In this subsection, we prove Theorem 2.4.
Because
\begin{align}
\begin{split}
h_{n,\ell}(y)&=yF\left(n+1,-n;\frac{1}{2};\frac{-\ell ^2}{4y^2-\ell ^2}\right)
=y\sum _{i=0}^n \frac{(n+1)_i(-n)_i}{(1/2)_i(1)_i}\left(\frac{-\ell ^2}{4y^2 -\ell ^2}\right)^i\\
&=\sum _{i=0}^n\frac{(n+1)_i }{(1/2)_i}\binom {n}{i}\ell ^{2i}\frac{y}{(4y^2 -\ell ^2)^i},
\end{split}
\end{align}
we have
\begin{align*}
\frac{d^2h_{n,\ell}}{dy^2}(y)&=
\sum _{i=0}^n\frac{(n+1)_i }{(1/2)_i}\binom {n}{i}\ell ^{2i}\frac{d^2}{dy^2}\frac{y}{(4y^2 -\ell ^2)^i}\\
&=\sum _{i=0}^n\frac{(n+1)_i }{(1/2)_i}\binom {n}{i}\ell ^{2i}
\frac{8iy\left((8i-4)y^2+3\ell ^2\right)}{(4y^2 -\ell ^2)^{i+2}}.
\end{align*}
Therefore, $h_{n,\ell}(y)$ is positive and convex on $\{y\in \R _{>0};y>\ell /2\}$, and
$$
\lim _{y\rightarrow \ell /2}h_{n,\ell}(y)=\lim _{y\rightarrow \infty}h_{n,\ell}(y)=\infty.
$$
Hence, there is a unique zero point $\mu _{n,\ell}$ of $h'_{n,\ell}(y)$.
Moreover, if we put $\nu _{n,\ell}:=h_{n,\ell}(\mu _{n,\ell})$,
then $\nu _{n,\ell}$ is the minimum of $h_{n,\ell}(y)$.

The role of $\nu _{n,\ell}$ and $minimality$ of
$D_{2m+1,\ell}(a,y_{m-1},y_{m-2},$ $\cdots ,y_0,y_0,\cdots ,y_{m-2},$ $ y_{m-1}, a)$
with $y_n=h_{n,\ell}(y_{a,m,\ell}^{\pm})$ for $n=0, \cdots ,m-1$
are obtained similarly as in the case $m=1,2,3$.

\begin{rmk}
We saw that
$h_{n,\ell}(y)$ is positive and convex on $\{y\in \R _{>0};y>\ell /2\}$ in the above.
We consider other properties of $h_{n,\ell}(y)$.

From (3.15), we have $h_{n+1,\ell}(y)>h_{n,\ell}(y)$ for $n\in {\Bbb N}$.
In addition, because
\begin{align*}
\frac{dh_{n,\ell}}{dy}(y)&=
\sum _{i=0}^n\frac{(n+1)_i }{(1/2)_i}\binom {n}{i}\ell ^{2i}\frac{d}{dy}\frac{y}{(4y^2 -\ell ^2)^i}\\
&=\sum _{i=0}^n\frac{(n+1)_i }{(1/2)_i}\binom {n}{i}\ell ^{2i}\frac{(4-8i)y^2-\ell ^2}{(4y^2 -\ell ^2)^{i+1}},
\end{align*}
we have $h' _{n,\ell}(y)>h' _{n+1,\ell}(y)$ for $n\in {\Bbb N}$.
Therefore, we obtain 
\begin{lmm}
$\mu _{n+1,\ell}> \mu _{n,\ell}$ and $\nu _{n+1,\ell}> \nu _{n,\ell}$ for $n\in {\Bbb N}$.
In particular, if $h' _{m,\ell}(y) >0 $ at a point $y$, then $h' _{n,\ell}(y) >0 $ for $n=0,1,\cdots ,m-1$.
\end{lmm}
\end{rmk}
\section{A proof of Theorem 2.5}
In this section, we prove Theorem 2.5.
For this, we investigate whether
$D_{2m+1,\ell}(a,$ $y_{m-1},$ $y_{m-2},$ $\cdots ,y_0,y_0,\cdots ,y_{m-2}, y_{m-1}, a)$
with $y_n=h_{n,\ell}(y_{a,m,\ell}^{\pm})$ for $n=0, \cdots ,m-1$\\ are stable.

We denote $S_{\ell}(s,t)$ by $S(s,t)$ for simplicity.
\subsection{Elements of the Hessian matrix}
We give an expression of each element of the Hessian matrix of 
the function $T_{a,m,\ell}(y_0,$ $ y_1, \cdots ,y_{m-1})$ at $(y,h_{1,\ell}(y),\cdots ,h_{m-1,\ell}(y))$
with $a=h_{m,\ell}(y)$.
We call this Hessian matrix $H_m(y)=(H_{i,j}(y))_{i, j=1,2,\cdots ,m}$.
Recalling that 
$$
T_{a,m,\ell}(y_0, y_1, \cdots ,y_{m-1})=y_0\ell+S(y_0,y_1)+S(y_1,y_2)+\cdots +S(y_{m-1},a),
$$
we see that the elements of $H_m(y)$ are expressed as
\begin{align*}
H_{1,1}(y)&=\frac{\partial ^2 S}{\partial s^2}(h_{0,\ell}(y),h_{1,\ell}(y)),\\
H_{i,i}(y)&=\frac{\partial ^2 S}{\partial t^2}(h_{i-2,\ell}(y),h_{i-1,\ell}(y))
+\frac{\partial ^2 S}{\partial s^2}(h_{i-1,\ell}(y),h_{i,\ell}(y))
\end{align*}
for $i=2,3,\cdots , m-1$,
$$
H_{m,m}(y)=\frac{\partial ^2 S}{\partial t^2}(h_{m-2,\ell}(y),h_{m-1,\ell}(y))
+\frac{\partial ^2 S}{\partial s^2}(h_{m-1,\ell}(y),a),
$$
$$
H_{i,i+1}(y)=H_{i+1,i}(y)=\frac{\partial ^2 S}{\partial s \partial t}(h_{i-1,\ell}(y),h_{i,\ell}(y))
$$
for $i=1, 2, \cdots , m-1$, and
$$
H_{i,j}(y)\equiv 0
$$
if $|i-j|\geq 2$.

Now, we consider the following function:
$$
\widehat{T}_{m,\ell}(y_0, y_1, \cdots ,y_{m-1},y_{m})=y_0\ell+S(y_0,y_1)+S(y_1,y_2)+\cdots +S(y_{m-1},y_{m})
$$
and define 
$\widehat{H}_m(y)=(\widehat{H}_{i,j}(y))_{i, j=1,2,\cdots ,m}$ as
$$
\widehat{H}_{i,j}(y):=\frac{\partial ^2 \widehat{T}_{m,\ell}}{\partial y_{i-1} \partial y_{j-1}}
(y, h_{1,\ell}(y), \cdots , h_{m-1,\ell}(y), h_{m,\ell}(y))
$$
for $i, j=1,2, \cdots m$.
Then, we see $\widehat{H}_{i,j}(y)=H_{i,j}(y)$ except for $(i,j)=(m,m)$ and
$$
\widehat{H}_{m,m}(y)=\frac{\partial ^2 S}{\partial t^2}(h_{m-2,\ell}(y),h_{m-1,\ell}(y))
+\frac{\partial ^2 S}{\partial s^2}(h_{m-1,\ell}(y),h_{m,\ell}(y)).
$$
Moreover, we notice that if $y=y_{a,m,\ell}^{+}$ or $y=y_{a,m,\ell}^{-}$, then $H_{m,m}(y)=\widehat{H}_{m,m}(y)$.
Therefore, for investing  whether
$D_{2m+1,\ell}(a,y_{m-1},y_{m-2},$ $\cdots ,y_0,y_0,\cdots ,y_{m-2}, y_{m-1}, a)$
with $y_n=h_{n,\ell}(y_{a,m,\ell}^{\pm})$ for $n=0, \cdots ,m-1$ are stable,
we evaluate $\widehat{H}_{m}(y)$ instead of $H_{m}(y)$.
\begin{rmk}
From the construction of $h_{m,\ell}(y)$, we see
\begin{gather}
\frac{\partial \widehat{T}_{m, \ell}}{\partial y_i}(y, h_{1, \ell}(y), \cdots h_{m-1, \ell}(y), h_{m, \ell}(y))
\equiv 0
\end{gather}
for $0\leq i \leq m-1$.
\end{rmk}
\subsection{An expression of the determinant of $\widehat{H}_{m}(y)$}
We consider the determinant of $\widehat{H}_{m}(y)$ for $m\geq 1$.
Note that each $\det\widehat{H}_{m}(y)$ satisfies the following recurrence relation:
\begin{gather}
\det \widehat{H}_{m}(y)=\widehat{H}_{m,m}(y)\det\widehat{H}_{m-1}(y)
-\left(\widehat{H}_{m-1,m}(y)\right)^2 \det\widehat{H}_{m-2}(y)
\end{gather} 
for $m\geq 3$. In this subsection, we prove 
\begin{lmm}
$$
\det \widehat{H}_{m}(y)=
(-1)^m \widehat{H}_{1,2}(y)\widehat{H}_{2,3}(y)\cdots \widehat{H}_{m,m+1}(y)h^{\prime}_{m,\ell}(y)
$$
for $m\geq 1$.
We remark that $\widehat{H}_{m}(y)$ does not have $\widehat{H}_{m,m+1}(y)$ as its element.
Here, we take $\widehat{H}_{m,m+1}(y)$ as
\begin{gather*}
\frac{\partial ^2 \widehat{T}_{m+1,\ell}}
{\partial y_{m-1}y_{m}}(y, h_{1,\ell}(y),\cdots ,h_{m-1,\ell}(y), h_{m,\ell}(y), h_{m+1, \ell}(y)).
\end{gather*}
\end{lmm}
We prove this lemma by induction.

First, we consider in the case that $m=1$. Because
\begin{align*}
\frac{\partial \widehat{T}_{1,\ell}}{\partial y_0}(y, h_{1,\ell}(y))\equiv 0
\end{align*}
from (4.1), we have
\begin{gather*}
0\equiv \frac{d}{dy}\left(\frac{\partial \widehat{T}_{1,\ell}}{\partial y_0}(y, h_{1,\ell}(y))\right)
=\frac{\partial ^2 \widehat{T}_{1,\ell}}{\partial y_0^2}(y, h_{1,\ell}(y))
+\frac{\partial ^2 \widehat{T}_{1,\ell}}{\partial y_0y_1}(y, h_{1,\ell}(y))\cdot h^{\prime}_{1,\ell}(y).
\end{gather*}
Observing that
\begin{gather*}
\frac{\partial ^2 \widehat{T}_{1,\ell}}{\partial y_0y_1}(y, h_{1,\ell}(y))=
\frac{\partial ^2 \widehat{T}_{2,\ell}}{\partial y_0y_1}(y, h_{1,\ell}(y), h_{2,\ell}(y))=
\widehat{H}_{1,2}(y),
\end{gather*}
we obtain Lemma 4.2 in the case that $m=1$.

Second, we consider in the case that $m=2$.
Since 
\begin{align*}
\frac{\partial \widehat{T}_{2,\ell}}{\partial y_1}(y, h_{1,\ell}(y), h_{2, \ell}(y))\equiv 0,
\end{align*}
we have
\begin{gather*}
0\equiv \frac{d}{dy}\left(\frac{\partial \widehat{T}_{2,\ell}}{\partial y_1}(y, h_{1,\ell}(y), h_{2,\ell}(y))\right)
=\widehat{H}_{1,2}(y)+\widehat{H}_{2,2}(y)\cdot h^{\prime}_{1,\ell}(y)+
\widehat{H}_{2,3}(y)\cdot h^{\prime}_{2,\ell}(y).
\end{gather*}
and consequently
\begin{align*}
(-1)^2\widehat{H}_{1,2}(y)\widehat{H}_{2,3}(y)\cdot h^{\prime}_{2,\ell}(y)
&=-\left(\widehat{H}_{1,2}(y)\right)^2-\widehat{H}_{1,2}(y)\widehat{H}_{2,2}(y)\cdot h^{\prime}_{1,\ell}(y)\\
&=-\left(\widehat{H}_{1,2}(y)\right)^2+\widehat{H}_{2,2}(y)\det \widehat{H}_1(y)=\det \widehat{H}_2(y)
\end{align*}
from the case $m=1$. 

Finally, we assume that this lemma holds for $1,2, \cdots m-1$.
We notice that
\begin{align*}
0&\equiv
\frac{d}{dy}\left(\frac{\partial \widehat{T}_{m,\ell}}{\partial y_{m-1}}(y, h_{1,\ell}(y),
\cdots h_{m-1,\ell}(y), h_{m,\ell}(y))\right)\\
&=\widehat{H}_{m-1,m}(y)\cdot h^{\prime}_{m-2,\ell}(y)+\widehat{H}_{m,m}(y)\cdot h^{\prime}_{m-1,\ell}(y)+
\widehat{H}_{m,m+1}(y)\cdot h^{\prime}_{m,\ell}(y).
\end{align*}
Therefore, by assumption and (4.2), we obtain 
\begin{flalign*}
&(-1)^m \widehat{H}_{1,2}(y)\cdots 
\widehat{H}_{m-1,m}(y) \widehat{H}_{m,m+1}(y)\cdot h^{\prime}_{m,\ell}(y)\\
&=(-1)^{m-1}\widehat{H}_{1,2}(y)\cdots \widehat{H}_{m-1,m}(y)\widehat{H}_{m,m}(y)\cdot h^{\prime}_{m-1,\ell}(y)\\
&\qquad +(-1)^{m-1}\widehat{H}_{1,2}(y)\cdots 
\widehat{H}_{m-1,m}(y)\widehat{H}_{m-1,m}(y)\cdot h^{\prime}_{m-2,\ell}(y)\\
&=\widehat{H}_{m,m}(y)\det \widehat{H}_{m-1}(y)-\left(\widehat{H}_{m-1,m}(y) \right)^2 \det\widehat{H}_{m-2}(y)
=\det\widehat{H}_{m}(y).
\end{flalign*}
Thus, we could show Lemma 4.2.
\subsection{A proof of Theorem 2.5}
In this subsection, we prove Theorem 2.5. \\
Recall that $S(s,t)=(s+t)\sqrt{(t-s)^2+\ell ^2}$ and 
\begin{align*}
\frac{\partial ^2 S}{\partial s\partial t}(s,t)&=\frac{-\ell ^2 (s+t)}{\left((t-s)^2 +\ell ^2\right)^{3/2}},\\
\widehat{H}_{i,i+1}(y)=H_{i,i+1}(y)&=\frac{\partial ^2 S}{\partial s \partial t}(h_{i-1,\ell}(y),h_{i,\ell}(y))
\end{align*}
for $i=1, \cdots ,m$. Therefore, every $\widehat{H}_{i,i+1}(y)$ is negative on $\{y\in \R_{>0}; y>\ell /2\}$
and thus, the sign of $\det \widehat{H}_{m}(y)$ is equal to the sign of $h^{\prime}_{m,\ell}(y)$ by Lemma 4.2.
The following lemma is well known:
\begin{lmm}
A symmetric $m\times m$ matrix $A=(A_{i,j})_{i,j=1,2,\cdots ,m}$ is positive definite 
if and only if $\det A^{(n)} >0$ for any $n=1,2, \cdots ,m$, where $A^{(n)}:=(A_{i,j})_{i,j=1,2,\cdots ,n}$.
\end{lmm}
Note that $\widehat{H}_m^{(n)}(y)$ is equal to $\widehat{H}_n(y)$.
If $y>\mu _{m,\ell}$, i.e., $y=y_{a,m,\ell}^{+}$,
then $\det H_{m}(y)=\det \widehat{H}_{m}(y)>0$
and $\det H_{n}(y)=\det \widehat{H}_n(y)>0$ for $n=1,2,\cdots ,m-1$ from Lemma 3.2.
Hence, from Lemma 4.3, we see that $H_{m}(y_{a,m,\ell}^{+}) $ is positive definite and
$D_{2m+1,\ell}(a,y_{m-1},y_{m-2},$ $\cdots ,y_0,y_0,\cdots ,y_{m-2}, y_{m-1}, a)$
with $y_n=h_{n,\ell}(y_{a,m,\ell}^{+})$ for $n=0, \cdots ,m-1$ are stable
for $a>\nu _{m, \ell}$.

\section{Approximations of catenaries}
In this section, we see
that the polylines (2.2) considerably approximate catenaries,
by calculating some cases numerically.

Recall a catenary $C_{c}(t):=c\cosh (t/c)$.
The function $c\mapsto c\cosh(1/c)$ is positive, convex on $\R _{>0}$ 
and takes the unique minimum $\eta _{\infty}:=1.5088 \cdots$ at $\xi _{\infty}:=0.8335\cdots$.
Thus, if $a>\eta _{\infty}$, 
then there are two positive numbers $c_a^{\pm}$ with 
$c_a^{-}< \xi _{\infty} <c_{a}^+$ such that $c_a^{\pm}\cosh (1/c_{a}^{\pm})=a$.

Before approximating catenaries,
we investigate a relation between $\eta _{\infty}$ and \\$\nu _{m,2/(2m+1)}$ in the following subsection.

\subsection{A relation between $\eta _{\infty}$ and $\nu _{m,2/(2m+1)}$}
In this subsection, we show the following lemma:
\begin{lmm}
$\eta _{\infty}> \nu _{m,2/(2m+1)}$ for $m \geq 1$. 
\end{lmm}
It suffices to prove that 
\begin{gather}
y\cosh (1/y)-h_{m,2/(2m+1)}(y)>0
\end{gather}
for $y\geq 1/\sqrt{2}$ because the first derivative of $y\cosh (1/y)$ at $y=1/\sqrt{2}$ is less than zero.
The left hand side of the Formula (5.1) is equal to
\begin{flalign*}
&y\cosh \left(\frac{1}{y}\right)-yF\left( m+1,-m;\frac{1}{2};\frac{-1}{(2m+1)^2y^2-1}\right) \\
&=\sum _{i=0}^{\infty}\frac{1}{(1)_{2i}\, y^{2i-1}}
-\sum _{i=0}^{\infty}\frac{(m+1)_i}{(1/2)_i}\binom{m}{i}\frac{y}{((2m+1)^2y^2-1)^i}\\
&=\sum _{i=0}^{\infty}\left(
\frac{1}{(1)_{2i}\, y^{2i-1}}-\frac{(m+1)_i}{(1/2)_i}\binom{m}{i}\frac{y}{((2m+1)^2y^2-1)^i}
\right)\\
&=:\sum _{i=0}^{\infty}a(m,i,y).
\end{flalign*} 
Note that $a(m,0,y)=0$. Now, we investigate $a(m,i,y)$.
First, we evaluate $a(m,1,y)$ $+a(m,2,y)$. This is equal to
\begin{gather*}
\frac{8m(m+1)y^2(2y^2-1)(3y^2+1)+(y-1)^2(y+1)^2(12y^2+1)}{24y^3((2m+1)^2y^2-1)^2}.
\end{gather*}
Therefore, if $y\geq 1/\sqrt{2}$, then $a(m,1,y)+a(m,2,y)$ is positive. 
Next, we evaluate $a(m,i,y)$ for $i\geq 3$.
Because
$$
\frac{1}{(2m+1)^2y^2-1}> \frac{1}{4(m+2)(m-1)y^2}
$$
for $y>1/3$,
\begin{align*}
a(m,i,y)&>\frac{1}{(1)_{2i}\, y^{2i-1}}-\frac{(m+1)_i}{(1/2)_i}\binom{m}{i}\frac{y}{4^i(m+2)^i(m-1)^iy^{2i}}\\
&=\frac{1}{(1)_{2i}\, y^{2i-1}}\left(1-\frac{(m+i)!}{(m-i)!(m+2)^i(m-1)^i}\right).
\end{align*}
for $y\geq 1/\sqrt{2}$.
We show 
$$
1> \frac{(m+i)!}{(m-i)!(m+2)^i(m-1)^i}=:b(m,i)
$$
for $i\geq 3$.
We easily find that $b(m,i)>b(m,i+1)$ for $i\geq 3$.
Therefore, we have only to show $1> b(m,3)$.
This is shown by
\begin{align*}
b(m,3)&=\frac{(m+3)!}{(m-3)!(m+2)^3(m-1)^3}\\
&= \frac{m(m-2)\cdot (m+1)(m+3)}{(m+2)^2(m-1)^2}
<\frac{(m-1)^2\cdot (m+2)^2}{(m+2)^2(m-1)^2}=1.
\end{align*}
Hence, $a(m,i,y)$ is positive for $i\geq 3$ and $y\geq 1/\sqrt{2}$.
Collecting the above results, we see that (5.1) is correct for $y\geq 1/\sqrt{2}$.
\subsection{Approximations of catenaries}
By Lemma 5.1, if $a>\eta _{\infty}$,
there are $y_{a,m,2/(2m+1)^{\pm}}$ for $m\geq 1$,
that is, we can approximate $C_{c_a^{\pm}}(t)$
by the polylines (2.2).

For example, if $a=2$,
then 
\begin{align*}
\begin{split}
&y_{2,1,2/3}^+\approx 1.7338,\\
&y_{2,2,2/5}^+\approx 1.7101,\ h_{1,2/5}(y_{2,2,2/5}^+)\approx 1.8050,\\
&y_{2,3,2/7}^+\approx 1.7035,\ h_{1,2/7}(y_{2,3,2/7}^+)\approx 1.7518,\  h_{2,2/7}(y_{2,3,2/7}^+)\approx 1.8497,\cdots
\end{split}
\end{align*}
$c_2^-\approx 0.4701$ and
\begin{align*}
\begin{split}
&C_{c_2^+}(1/3)\approx 1.7295,\\
&C_{c_2^+}(1/5)\approx 1.7084,\ C_{c_2^+}(3/5)\approx 1.8038,\\
&C_{c_2^+}(1/7)\approx 1.7026,\ C_{c_2^+}(3/7)\approx 1.7510,\  C_{c_2^+}(5/7)\approx 1.8492, \cdots.
\end{split}
\end{align*}
Moreover,
\begin{align*}
\begin{split}
&y_{2,1,2/3}^-\approx 0.5150,\\
&y_{2,2,2/5}^-\approx 0.4856,\ h_{1,2/5}(y_{2,2,2/5}^-)\approx 0.8823,\\
&y_{2,3,2/7}^-\approx 0.4779,\ h_{1,2/7}(y_{2,3,2/7}^-)\approx 0.6655,\  h_{2,2/7}(y_{2,3,2/7}^-)\approx 1.1141,\cdots
\end{split}
\end{align*}
$c_2^+\approx 1.6966$ and
\begin{align*}
\begin{split}
&C_{c_2^-}(1/3)\approx 0.5933,\\
&C_{c_2^-}(1/5)\approx 0.5133,\ C_{c_2^-}(3/5)\approx 0.9078,\\
&C_{c_2^-}(1/7)\approx 0.4920,\ C_{c_2^-}(3/7)\approx 0.6794,\  C_{c_2^-}(5/7)\approx 1.1254, \cdots.
\end{split}
\end{align*}
So, we see that the polylines (2.2) form nearly equal catenaries $C_{c_a^{\pm}|_{(-1,1)}}$ 
for $m$ large enough.
We draw two polylines whose vertices are specified by the sequence
\begin{gather*}
\begin{split}
\left(
(-1, a), \left(\frac{1-2m}{2m+1}, h_{m-1,\frac{2}{2m+1}}\left( y_{a,m,\frac{2}{2m+1}}^{-} \right)\right) , \cdots ,
\left(\frac{-1}{2m+1}, y_{a,m,\frac{2}{2m+1}}^{-}\right) ,\right. \\
\left.
\left(\frac{1}{2m+1}, y_{a,m,\frac{2}{2m+1}}^{-}\right),\cdots, 
\left(\frac{2m-1}{2m+1}, h_{m-1,\frac{2}{2m+1}}\left( y_{a,m,\frac{2}{2m+1}}^{-} \right)\right),(1,a)
\right)
\end{split}
\end{gather*}
in cases that $(a,m)=(2,5)$ and $(2,10)$ below:
\begin{figure}[bpht]
 \begin{center}
  \includegraphics[width=100mm]{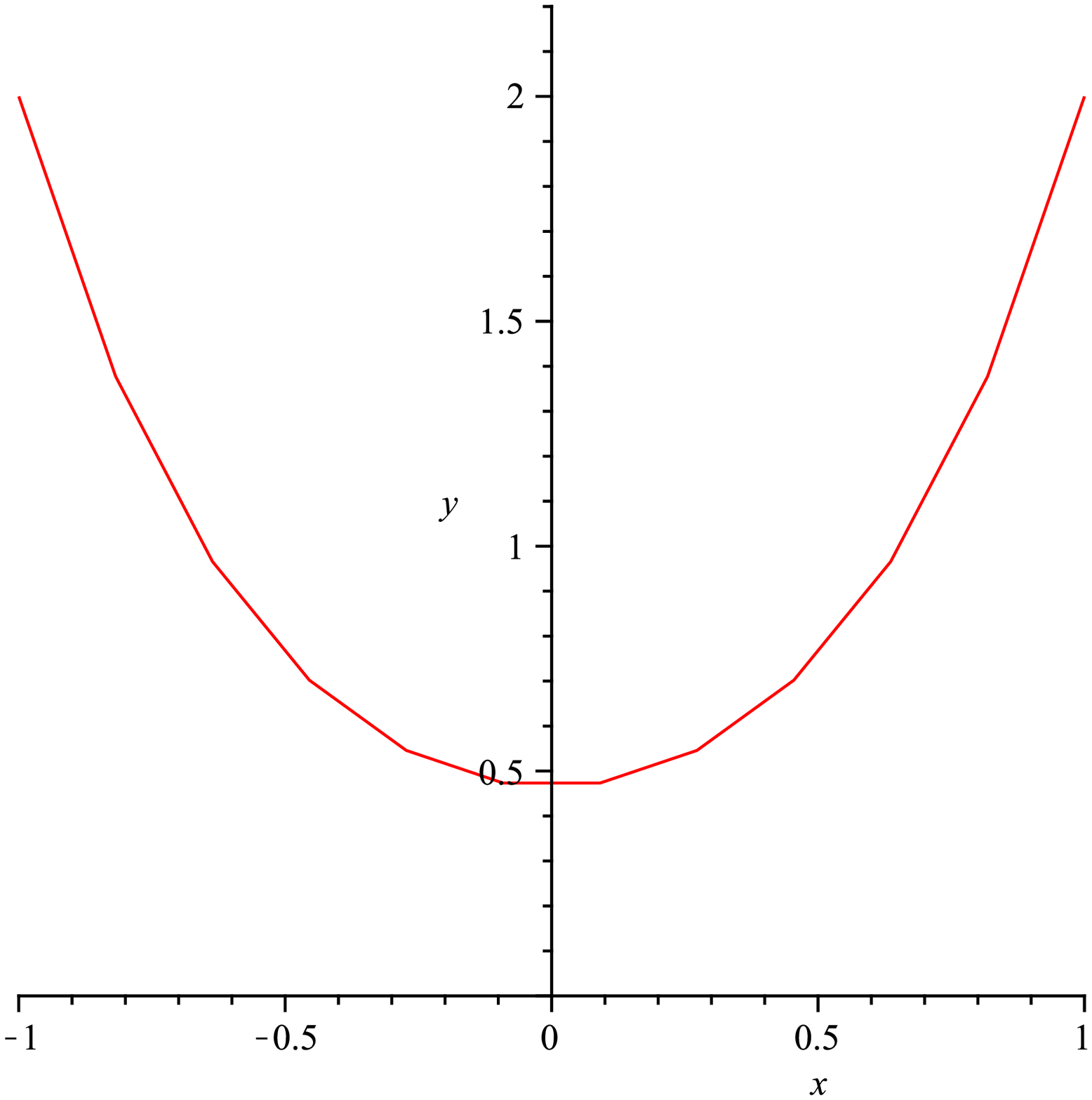}
 \end{center}
 \caption{$(a,m)=(2,5)$}
 \label{fig:one}
\end{figure}
\begin{figure}[bpht]
 \begin{center}
  \includegraphics[width=100mm]{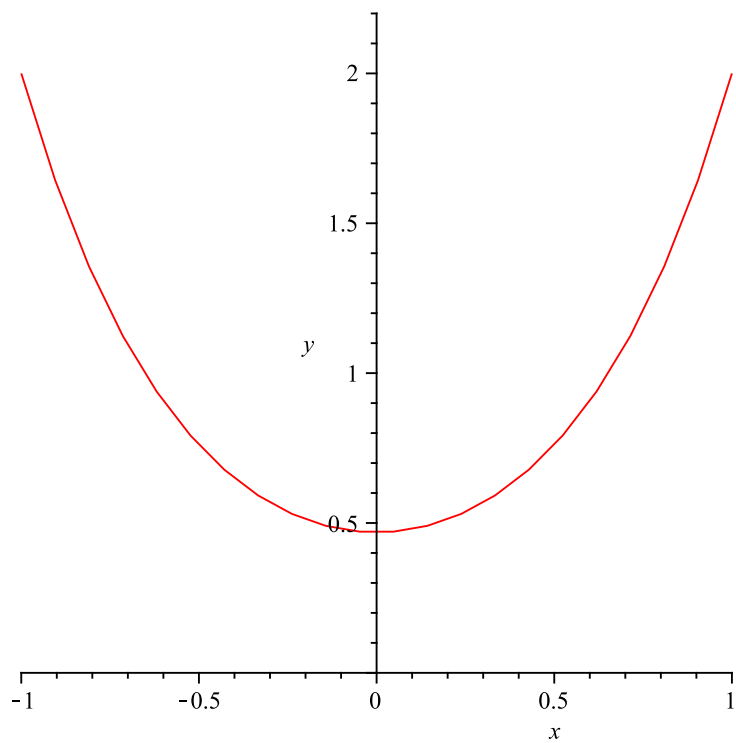}
 \end{center}
 \caption{$(a,m)=(2,10)$}
 \label{fig:two}
\end{figure}

\medskip
\begin{flushleft}
Akihito Ebisu\\
Department of Mathematics\\
Kyushu University\\
Nishi-ku, Fukuoka 819-0395\\
Japan\\
a-ebisu@math.kyushu-u.ac.jp
\end{flushleft}

\begin{flushleft}
Yoshiroh Machigashira\\
Division of Mathematical Sciences \\
Osaka Kyoiku  University\\
4-698-1, Asahigaoka, Kashiwara, Osaka 582-8582\\
Japan\\
machi@cc.osaka-kyoiku.ac.jp
\end{flushleft}

\end{document}